\documentclass[12pt,a4paper]{amsart}

\theoremstyle{plain}

\usepackage{enumerate, amssymb}

\advance\hoffset-5mm \advance\textwidth40mm
\input{diagrams.tex}
\diagramstyle[scriptlabels,height=8mm,width=8mm]

\def\bdi{\begin{diagram}}
\def\edi{\end{diagram}}

\theoremstyle{plain}

\newtheorem{thm}{Theorem}[section]
\newtheorem{cor}[thm]{Corollary}
\newtheorem{lem}[thm]{Lemma}
\newtheorem{prop}[thm]{Proposition}
\theoremstyle{definition}
\newtheorem{defi}[thm]{Definition}
\newtheorem{defis}[thm]{Definitions}
\newtheorem{conj}[thm]{Conjecture}
\newtheorem{conv}[thm]{Convention}
\newtheorem{nota}[thm]{Notation}
\newtheorem{rem}[thm]{Remark}
\newtheorem{rems}[thm]{Remarks}
\newtheorem{exa}[thm]{Example}
\newtheorem{exas}[thm]{Examples}
\newtheorem{prob}[thm]{Problem}
\newtheorem{probs}[thm]{Problems}
\newtheorem{ques}[thm]{Question}
\newtheorem{sit}[thm]{}

\newcommand{\IVF}{ \operatorname{{\rm IVF}}}
\newcommand{\AVF}{ \operatorname{{\rm AVF}}}

\newcommand{\Lie}{ \operatorname{{\rm Lie}}}
\newcommand{\LieA}{ \operatorname{{\rm Lie_{alg}}}}
\newcommand{\LieAO}{ \operatorname{{\rm Lie_{alg}^\omega}}}
\newcommand{\Ker}{ \operatorname{{\rm Ker}}}

\newcommand{\Aut}{ \operatorname{{\rm Aut}}}

\renewcommand{\epsilon}{\varepsilon}

\def\and{\quad\mbox{and}\quad}

\newcommand{\C}{\ensuremath{\mathbb{C}}}

\newcommand{\R}{\ensuremath{\mathbb{R}}}

\newcommand{\Z}{\ensuremath{\mathbb{Z}}}

\newcommand{\sgoth}{{\ensuremath{\mathfrak{s}}}}
\newcommand{\lgoth}{{\ensuremath{\mathfrak{l}}}}

\newcommand{\cB}{{\ensuremath{\mathcal{B}}}}
\newcommand{\cL}{{\ensuremath{\mathcal{L}}}}

\newcommand{\cS}{{\ensuremath{\mathcal{S}}}}

\newcommand{\cC}{{\ensuremath{\mathcal{C}}}}

\newcommand{\cN}{{\ensuremath{\mathcal{N}}}}

\newcommand{\cW}{{\ensuremath{\mathcal{W}}}}
\newcommand{\cZ}{{\ensuremath{\mathcal{Z}}}}

\renewcommand{\rho}{\varrho}

\def\bals#1\eals{\begin{align*}#1\end{align*}}
\def\bal#1\eal{\begin{align}#1\end{align}}

\def\SAut{\mathop{\rm SAut}}

\def\AA{{\mathbb A}}

\renewcommand{\phi}{\varphi}

\newcommand{\bnum}{\begin{enumerate}}
\newcommand{\enum}{\end{enumerate}}
\renewcommand{\emptyset}{\varnothing}

\newcommand{\brem}{\begin{rem}}
\newcommand{\brems}{\begin{rems}}
\newcommand{\erem}{\end{rem}}
\newcommand{\erems}{\end{rems}}
\newcommand{\bprob}{\begin{prob}}
\newcommand{\eprob}{\end{prob}}
\newcommand{\bprobs}{\begin{probs}}
\newcommand{\eprobs}{\end{probs}}
\newcommand{\bques}{\begin{ques}}
\newcommand{\eques}{\end{ques}}
\newcommand{\bexa}{\begin{exa}}
\newcommand{\bexas}{\begin{exas}}
\newcommand{\eexa}{\end{exa}}
\newcommand{\eexas}{\end{exas}}
\newcommand{\bdefi}{\begin{defi}}
\newcommand{\edefi}{\end{defi}}
\newcommand{\bdefis}{\begin{defis}}
\newcommand{\edefis}{\end{defis}}
\newcommand{\bcor}{\begin{cor}}
\newcommand{\ecor}{\end{cor}}
\newcommand{\blem}{\begin{lem}}
\newcommand{\elem}{\end{lem}}
\newcommand{\bconv}{\begin{conv}}
\newcommand{\econv}{\end{conv}}
\newcommand{\bconj}{\begin{conj}}
\newcommand{\econj}{\end{conj}}
\newcommand{\bprop}{\begin{prop}}
\newcommand{\eprop}{\end{prop}}
\newcommand{\bthm}{\begin{thm}}
\newcommand{\ethm}{\end{thm}}
\newcommand{\bnota}{\begin{nota}}
\newcommand{\enota}{\end{nota}}
\newcommand{\bsit}{\begin{sit}}
\newcommand{\esit}{\end{sit}}
\newcommand{\be}{\begin{equation}}
\newcommand{\ee}{\end{equation}}
\newcommand{\bproof}{\begin{proof}}
\newcommand{\eproof}{\end{proof}}
\def\ba{\begin{array}}
\def\ea{\end{array}}

\title[Algebraic (Volume) Density Property for Affine Homogeneous Spaces]
{Algebraic (Volume) Density Property for Affine Homogeneous Spaces}
\author{Shulim Kaliman}
\address{Department of Mathematics\\
University of Miami\\
Coral Gables, FL 33124 \ \ USA}
\email{kaliman@math.miami.edu}
\author{Frank Kutzschebauch}
\address{Mathematisches Institut \\Universit\"at Bern
 \\Sidlerstr. 5
 \\ CH-3012 Bern, Switzerland}
\email{Frank.Kutzschebauch@math.unibe.ch}

\thanks{{\bf Acknowledgements:} This research was started during a visit
 of the second author to the University of Miami,
Coral Gables. He thanks this institutions for its generous
support and excellent working conditions. The research of the
second author was  partially supported by Schweizerische
Nationalfonds grant No.  200021-153120/1}


\begin{document}
\begin{abstract} Let $X$ be   a connected affine homogenous space of a linear algebraic group $G$ over $\C$.
(1) If $X$ is different from a line or a torus we show that the space of
all algebraic vector fields on $X$  coincides with the Lie algebra
generated by complete algebraic vector fields on $X$.
(2) Suppose that $X$ has a $G$-invariant volume form $\omega$. We prove that the space of
all divergence-free (with respect to $\omega$) algebraic vector fields on $X$  coincides with the Lie algebra
generated by divergence-free complete algebraic vector fields on $X$ (including the cases when $X$ is a line or a torus). 
The proof of these results requires new criteria for algebraic (volume) density property based on so called module generating pairs.
 \end{abstract}
\maketitle \vfuzz=2pt

\vfuzz=2pt

\section{Introduction}

In this paper we continue to study algebraic (volume) density property in the framework of the Andersen-Lempert theory
for review of which we refer to \cite{Ro}, \cite{KaKu3}, and \cite{For}). 
Recall that a smooth complex affine algebraic variety $X$ possesses an algebraic volume form if it
admits a nowhere vanishing top algebraic differential form $\omega$ (or, equivalently, a nowhere vanishing section 
of its canonical bundle). We consider the set $\IVF (X)$ of complete
algebraic vector fields, its subset  $\IVF_\omega (X)$ that consists of fields of zero divergence with respect to $\omega$,
and the Lie algebra $\LieA (X)$ (resp.  $\LieAO (X)$) generated by $\IVF (X)$ (resp. $\IVF_\omega (X)$).

\bdefi\label{nc.08.20.10}
We say that
a smooth affine algebraic variety $X$ (over $\C$) has the algebraic  density property (ADP) if the Lie algebra $\LieA (X)$ generated by the set $\IVF (X)$ of complete
algebraic vector fields coincides with the space $\AVF (X)$ of all algebraic vector fields on $X$. Similarly in the presence of $\omega$ we can speak about 
the algebraic volume density property (AVDP) if $\LieAO (X)$ coincides with the space $\AVF_\omega (X)$ of all algebraic
vector fields that have $\omega$-divergence zero (and we call such fields divergence-free).
\edefi

In \cite{DDK} and \cite{KaKu} the following facts were established.

\bthm\label{int1}
Let $G$ be a connected linear algebraic group, $R$ be a proper reductive subgroup of $G$, and $X$ be  the homogeneous space $X=G/R$. 

{\rm (1)} Suppose that $X=G/R$ is not isomorphic to a line or a torus $(\C^*)^k$. Then $X$ has ADP.

{\rm (2)} Suppose that $X$
has a left-invariant (with respect to the natural $G$-action) volume form $\omega$.
Then $X$ has AVDP with respect to $\omega$.
\ethm

The requirement that $R$ is reductive enables us to insure that $G/R$ is affine by the Matsushima theorem \cite{Ma}.
However not every affine homogeneous space of a linear algebraic group can be presented in this form (e.g., see \cite{Snow}).
Therefore the aim of this paper is to strengthen Theorem \ref{int1} as follows. 

\bthm\label{main'}
Let  $X$ be  a connected affine homogeneous space of a linear algebraic group $G$.

{\rm (1)} Suppose that $X$ is not isomorphic to a line or a torus $(\C^*)^k$. Then $X$ has ADP.

{\rm (2)} Suppose that $X$
has a left-invariant (with respect to the natural $G$-action) volume form $\omega$.
Then $X$ has AVDP with respect to $\omega$.

\ethm

The proof of this theorem requires new methods for establishing ADP and AVDP. It is the special case when the ``semisimple part" of the homogeneous space in question is
$SL_2 / \C^*$ which we are unable to attack with methods developed earlier. If  ``semisimple part" is different from $SL_2 / \C^*$  (so-called general case)
then the old technique works. 

The paper is organized as follows. In section \ref{cri} we recall the main
criterion   (Theorem \ref{mcd1}) for proving AVDP based on the concept of semi-compatible pairs developed in \cite{KaKu}. Then we introduce a more general notion of
a module generating pair and explain that Theorem \ref{mcd1} is true for such pairs as well. In section \ref{sc} we recall known facts about semi-compatible pairs and an important existence result for them from \cite{DDK} which will be used in the general case. In section \ref{has} we recall results 
about the structure of  affine homogeneous spaces from the work of Snow \cite{Snow}. In the next section we establish condition (A) of Theorem \ref{mcd1} for the general case and  in the section thereafter condition (B)  
for all affine homogeneous spaces. We summarize these results in the next section, proving part (2) of our main theorem in the general case. 
Section \ref{newtech1} gives a new technique of proving ADP.
It is a generalization of the technique of compatible pairs developed in \cite{KaKu1} and it
 is  well adapted  to prove ADP for the total space of  certain Zariski locally trivial fibrations with base and fibre being flexible.
In section \ref{newtech2} we present a criterion for proving the existence
of a module generating pair in the case of similar fibrations. 
In the next section we can complete part (2) of our main theorem by closely analyzing the special case and pointing out how the techniques in sections \ref{newtech1} and \ref{newtech2}
can be applied. The last section is devoted to the proof of part (1) of our main theorem which in the general case uses methods developed in our earlier work and in the 
special case again appeals to our new techniques from section \ref{newtech1}.

\section{Main Criterion for AVDP}\label{cri}

The  aim of this section is to recall the main criterion for AVDP.

\bnota\label{mcn1}
 For the rest of the paper $X$ is always a  {\bf smooth affine irreducible algebraic
variety over $\C$}.
We consider often a situation when $X$ is equipped also with  an  {\bf effective}  algebraic action of a group $\Gamma$. We call such an $X$ a  $\Gamma$-variety
and for the rest of the paper $\Gamma$ is always  {\bf finite}.
The ring of regular $\Gamma$-invariant functions will be denoted by $\C [X, \Gamma ]$; it is naturally isomorphic to the ring
$\C [X/\Gamma ]$  of regular functions on the quotient space  $X/\Gamma$. 

\enota

\bdefi\label{mcd1} 
(1) Recall that a holomorphic vector field $\xi$ on $X$ is called complete
if there is a holomorphic $\C_+$-action $\Phi : \C \times X \to X$
such that $\xi (f)= {\frac{\rm d} {{\rm d} t}} f \circ \Phi (t, * )|_{t=0}$ for every  holomorphic function $f$ on $X$.
This action $\Phi$ is called the phase flow of $\xi$. When $\Phi$ is an algebraic $\C_+$-action
the field $\xi$ is called locally nilpotent, and  when $\Phi$ factors through $\C^* \times X$ and generates
an algebraic $\C^*$-action then $\xi$ is called semi-simple.

(2)  Let $\xi$  and $\eta$ be
nontrivial complete algebraic vector fields on a    $\Gamma$-variety $X$ which are $\Gamma$-invariant.
We say that the pair $(\xi , \eta )$ is $\Gamma$-semi-compatible if\\
\begin{center} {\em the span of
$ (\Ker \xi \cap \C (X, \Gamma ))\cdot (\Ker \eta \cap \C (X, \Gamma )) $  contains a nonzero ideal of  $\C [X, \Gamma ]$.}\\[2ex]
\end{center}
 The largest ideal contained in the span will be called the associate $\Gamma$-ideal of the pair $(\xi , \eta )$.
In the case of a trivial $\Gamma$ we say that the pair $(\xi , \eta)$ is semi-compatible. In this terminology $\Gamma$-semi-compatibility
of $(\xi , \eta )$  is equivalent to  semi-compatibility of $(\xi' , \eta')$ where $\xi'$ and $\eta'$ are the vector fields on $X/\Gamma$
induced  by $\xi$ and $\eta$.

(3) A $\Gamma$-semi-compatible pair $( \xi , \eta )$ is called $\Gamma$-compatible, if $\xi$ is locally nilpotent and there exists $a \in \Ker \eta \cap \C [X, \Gamma ]$
such that $\xi (a) \ne 0$ and $\xi^2 (a)=0$ (i.e. the $\xi$-degree of $a$ is 1). In the case of trivial $\Gamma$-action we say that the pair is compatible.

\edefi

\bexa\label{mce1} (1)  Consider the direct product $X=X_1\times X_2$ of affine algebraic varieties and let $\xi_i$ be a nontrivial complete algebraic
vector field on $X_i$. Suppose that $\eta_i$ is the natural lift of $\xi_i$ to $X$ and treat $\C [X_i]$ as a subring of  $\C [X]$.
Note that $\eta_1$ and $\eta_2$ are semi-compatible since $\C [X_2] \subset \Ker \eta_1$ and  $\C [X_1] \subset \Ker \eta_2$
and thus $\C [X_1] \cdot \C [X_2] \subset \Ker \eta_1 \cdot \Ker \eta_2$.

(2) If in (1) $\xi_1$ is locally nilpotent then the pair is automatically compatible. Indeed, there is always $a \in \C [X_1]$
of $\xi$-degree 1 and it belongs to the kernel of $\xi_2$.

\eexa

\bnota\label{mcn2}
(1) Suppose that $X$ is equipped with an algebraic volume form $\omega$.
Let $\cC_{k} (X)$ be the space of
algebraic differential $k$-forms on $X$ and $\cZ_{k}(X)$ and $\cB_{k} (X)$ be its subspaces of
closed and exact $k$-forms respectively. If $\dim X =n$ then there exists an isomorphism $\Theta : \AVF_\omega (X) \to \cZ_{n-1}(X)$
given by the formula $\xi \to  \iota_\xi \omega$ where $ \iota_\xi \omega$ is the interior product of $\omega$ and $\xi \in \AVF_\omega (X)$.

(2) In addition to (1) suppose that $X$ be a
$\Gamma$-variety.
Then for every $\gamma \in \Gamma$ one has
$$\omega \circ \gamma = \chi_\omega (\gamma ) \omega$$ where
$\chi_\omega$  is a map from  $\Gamma$  into the group of invertible regular functions on $X$.
In general $\chi_\omega (\gamma )$ may not be constant but we shall show 
later that for our proof it is enough to consider  the case in which $\chi_\omega$ is a character of $\Gamma$.
We denote by $\LieAO (X , \Gamma)$
the Lie algebra generated by $\Gamma$-invariant complete algebraic vector fields of $\omega$-divergence zero.
Similarly $\AVF_\omega (X, \Gamma )$ will be the space of algebraic $\Gamma$-invariant divergence-free vector fields.
Such fields can be also viewed as elements of $\AVF_\omega (X/\Gamma )$. By  $\cC_k (X, \Gamma )$  we denote the space
of algebraic $k$-forms  $\alpha$ such that
$$\alpha \circ \gamma = \chi_\omega (\gamma ) \alpha $$ for every $\gamma \in \Gamma$. 
By   $\cZ_k (X, \Gamma )$
and $\cB_k (X, \Gamma )$ we denote
the subspaces of closed and exact forms in $\cC_k (X, \Gamma )$.
The restriction of $\Theta : \AVF_\omega (X) \to \cZ_{n-1}(X)$ generates
an isomorphism between $\AVF_\omega (X, \Gamma )$ and  $\cZ_{n-1} (X, \Gamma )$.

(3) The set of complete $\Gamma$-invariant  divergence-free algebraic vector fields will
be denoted by $\IVF_\omega (X, \Gamma)$ and the Lie algebra generated by it by $\LieAO (X, \Gamma)$.
 \enota




The following fact is our main tool for establishing AVDP.

\bthm\label{mct1} {\rm (\cite[Theorem 2]{KaKu})}  
 Let  $X$ be a $\Gamma$-variety equipped with an algebraic volume form $\omega$, $\chi_\omega$ be a character of $\Gamma$,
and let $(\xi_j , \eta_j )_{j=1}^k$ be pairs of
$\Gamma$-semi-compatible  divergence-free  vector fields.
Let $I_j$ be the $\Gamma$-ideal associated with $(\xi_j , \eta_j )$
and $I_j(x)= \{ f (x) | f \in I_j \}$ for $x \in X$. Suppose that\\

\noindent {\rm (A)}  \hspace{0.3cm} for every $x \in X$ the set $\{ I_j (x) \xi_j (x) \wedge \eta_j (x) \}_{j=1}^k$ generates  $\Lambda^2 T_{x}X$.  \\

Suppose also that the following condition is true\\
 
\noindent {\rm (B)} \hspace{0.3cm} the image of  $\Theta (\LieAO (X,\Gamma))$  under
De Rham homomorphism $\Phi_{n-1} : \cZ_{n-1} (X) \to H^{n-1} (X, \C)$ coincides with the subspace  $\Phi_{n-1}( \cZ_{n-1} (X, \Gamma))$  of $H^{n-1}(X , \C)$.\\

Then
$\Theta (\LieAO (X, \Gamma )) = \cZ_{n-1} (X, \Gamma)$  and therefore
$\LieAO (X, \Gamma ) = \AVF_\omega (X, \Gamma)$.

\ethm

\bdefi\label{mcd2} Following \cite[Definition 3.14]{KaKu} we say that in the case of equality $\LieAO (X, \Gamma ) = \AVF_\omega (X, \Gamma)$
the variety $X$ possesses $\Gamma$-AVDP (with respect to $\omega$).

\edefi 

\brem\label{mcr1} It is worth explaining the reason for Condition (A) in Theorem \ref{mct1}. Consider the differential map ${\rm d} : \cC_{n-2}(X, \Gamma ) \to \cB_{n-1}(X,\Gamma )$ and
let $(\xi_i , \eta_i )$ be a pair of
$\Gamma$-semi-compatible  divergence-free  vector fields associated with  $\Gamma$-ideal $I_i$. Then it was proven in \cite{KaKu} that
${\rm d}^{-1} (\Theta (\LieAO (X, \Gamma ))$ contains the $\C [X, \Gamma ]$-submodule of $\cC_{n-2}(X, \Gamma )$
that consists of all $(n-2)$-forms that can be presented as $f \iota_{\xi_i} \iota_{\eta_i} \omega$ with $f \in I_i$. Condition (A) guarantees that
$$\sum_i I_i \iota_{\xi_i} \iota_{\eta_i} \omega = \cC_{n-2} (X, \Gamma )$$ and therefore $\Theta (\LieAO (X, \Gamma ))$ contains $\cB_{n-1}(X, \Gamma)$.
In combination with Condition (B) this yields the equality  $\Theta (\LieAO (X, \Gamma ))=\cZ_{n-1}(X,\Gamma )$ and thus  $\Gamma$-AVDP for $X$. 
\erem

\bdefi\label{mcd3}  Let  $(\xi , \eta )$ be a pair of
$\Gamma$-invariant  divergence-free algebraic vector fields on $X$.    We say that   this pair $(\xi , \eta )$
is $\Gamma$-module generating if ${\rm d}^{-1} (\Theta (\LieAO (X, \Gamma ))$ contains the $\C [X, \Gamma ]$-submodule of $\cC_{n-2}(X, \Gamma )$
that consists of all $(n-2)$-forms that can be presented as $f \iota_{\xi} \iota_{\eta} \omega$ with $f$ in a  $\Gamma$-ideal $J$. The largest such ideal $J$  will be called the ideal associated with
the pair, or shortly the associated ideal. 
\edefi 

\brem\label{mcr2} Remark \ref{mcr1}  shows that a pair of $\Gamma$-semi-compatible  divergence-free complete vector fields is  $\Gamma$-module generating. Moreover is is easily seen  that Theorem \ref{mct1} remains valid if instead of requiring that pairs 
$(\xi_i , \eta_i )$ are $\Gamma$-semi-compatible  divergence-free  vector fields one assume  that they are $\Gamma$-module generating.

\erem

 We introduce this new definition of module generating fields (which by the way need not to be complete) since
in the special low dimensional case mentioned in the introduction the semi-compatibility alone does not work.

\section{Semi-compatibility}\label{sc}

Finding semi-compatible pairs of vector fields on $X$ is crucial and hence we
are going to exploit criteria developed in our previous papers some of which are collected in 
the next Proposition (see \cite[Proposition 3.4]{KaKu1},  \cite[Theorem 12]{DDK}, and  \cite[Remark 2.5, Propositions 2.4, and 2.13]{KaKu}).

\bprop\label{scp1}
{\rm (1)} Let $\xi_1$ be a locally nilpotent vector field on $X$ with a finitely generated kernel and let $\xi_2$ be either a similar field or a semi-simple one.
Suppose that $\rho_i : X \to X_i:={\rm Spec} \, \Ker \xi_i$ is the natural morphism and $\rho = (\rho_1, \rho_2): X \to X_1 \times X_2$.
Then the pair $(\xi_1, \xi_2)$ is semi-compatible if and only if $Z:=\rho (X)$ is closed in $X_1 \times X_2$ and $\rho$ induces a
finite birational morphism $X \to Z$.

{\rm (2)} 
Let $X$ admit a fixed point free non-degenerate $SL_2$-action
(i.e. general $SL_2$-orbits are of dimension 3)
and let $U_1$ (resp. $U_2$) be the unipotent subgroup of upper (resp. lower) triangular matrices in $SL_2$.
Suppose that $\xi_1$ and $\xi_2$ are the locally nilpotent
vector  fields   associated with the induced actions of the $\C_+$-groups $U_1$ and $U_2$ on $X$. Then $(\xi_1 , \xi_2 )$ is a compatible 
(and in particular semi-compatible) pair.

{\rm (3)} Let $\xi_1$ and $\xi_2$ be two semi-simple vector fields on $X$
such that they commute, i.e. they induce a $\C^*\times \C^*$-action on $X$. Suppose that
the algebraic quotient morphism $\pi : X \to W :=X//(\C^* \times \C^*)$ is smooth and all orbits of the action are two-dimensional.
Then the pair $(\xi_1 , \xi_2 )$ is semi-compatible.

{\rm (4)} Let $X$ be a  $\Gamma$-variety, $\xi_1$ and $\xi_2$
be semi-compatible  $\Gamma$-invariant vector fields on $X$. Suppose that $\xi_1'$ and $\xi_2'$ are the induced vector fields on $X'=X/\Gamma$
and one of the following conditions holds:

{\rm (i)} $\xi_1$ and $\xi_2$ are locally nilpotent and they generate a Lie algebra isomorphic to
$\sgoth \lgoth_2$ that induces a fixed point free non-degenerate algebraic $SL_2$-action on $X$;

{\rm (ii)} $[\xi_1 , \xi_2 ]=0$ and $\xi_1$ is a locally nilpotent   vector field with
a finitely generated kernel   while $\xi_2$ is
either locally nilpotent (also with a finitely generated kernel) or semi-simple.

Then the pair $(\xi' , \eta' )$ is semi-compatible, i.e. $(\xi , \eta)$ is $\Gamma$-semi-compatible.

\eprop


\bdefi\label{scd1} (1) For any vector field $\sigma$ on $X$  a field $f \sigma$ with $f \in \Ker \sigma$ is
called a replica of $\sigma$. 

(2) Suppose that $\cN$ is a collection of locally nilpotent vector fields and let $F$ be the subgroup of $\Aut (X)$
generated by the elements of phase flows of fields from $\cN$ and their replicas.
We say that $X$ is $F$-flexible if $F$ acts transitively on $X$, or, equivalently, for every $x \in X$ the set
$\{ \sigma |_x \,  | \sigma \in \cN \}$ generates $T_xX$  (see \cite{AFKKZ} for
properties of such varieties).

In the case when $\cN$ is the set of all locally nilpotent vector fields we call $X$ flexible and denote $F$ by $\SAut (X)$.
\edefi

The next claim is a combination of Lemma 5.1 and Remark 5.2 from \cite{KaKu}.

\bprop\label{flexible} Let $X$ be a flexible affine algebraic manifold that possesses a pair of semi-compatible vector fields.
Then $X$ satisfies Condition ($A$) of Theorem \ref{mct1} (with trivial $\Gamma$). 
\eprop

We need an equivariant version of Proposition \ref{flexible} which requires the following technical fact.

\blem\label{Gam} Let $\Gamma$ be a finite group acting on an affine algebraic variety $X$ 
and $\sigma$ be a locally nilpotent $\Gamma$-invariant vector field on $X$.
Then for a general point $x \in X$ and  every vectors $u$ and $v\in T_xX$ such that
$u,v$ and the value of the field $\sigma$ at $x$ are linearly independent, one has a regular $\Gamma$-invariant function
$f$ on $X$ for which $f(x)=0$, ${\rm d}f(u)=0$ and $ {\rm d}f(v) \ne 0$.

\elem

\bproof  Since $x$ is a general point we can suppose
that the image $a$ of $x$ under the quotient morphism $\tau : X \to A ={\rm Spec} \Ker \sigma$ is a general point, i.e. 
there is a neighborhood of $a$ in $A$ which is a smooth affine algebraic variety such that the restriction of $\tau$ over this
neighborhood is a trivial $\AA^1$-fibration. By linear independence
the images $u'$ and $v'$ of $u$ and $v$ under ${\rm d} \tau$ are non-proportional vectors in $T_aA$. 
Because of smoothness of $A$ at $a$ one can choose a function $g$ on $A$ for which $g(a)=0$, ${\rm d}g(u' )=0$ and $ {\rm d}g(v' ) \ne 0$.
Furthermore, because of commutativity $\Gamma$ acts on $A$ and this function $g$ can be chosen
as a lift a function on $A/\Gamma$ (since $A/\Gamma$ is also smooth at the image of
the general point $a$). This implies that the lift $f$ of $g$ to $X$ is $\Gamma$-invariant which yields the desired function.

\eproof

\bprop\label{Gamma} Let $X$ be a $\Gamma$-variety and $\varphi : X \to T$
be a morphism.  Let $\cN$ and $F$ be as in Definition \ref{scd1} with all elements of $\cN$ being $\Gamma$-invariant. Suppose that 
the action of $F$ preserves a fiber $Z$ of $\varphi$  with $\dim Z \geq 3$.
Furthermore, let $Z$ be $F$-flexible. Suppose 
$u\wedge v$  is a nonzero element of $\Lambda^2T_xZ$ at some $x\in Z$.
Consider the subgroup $F_\Gamma$ of $F$ that commutes with $\Gamma$
and the isotropy
group $F_\Gamma^x$ of $x$ in $F_\Gamma$.
Then the orbit of $u\wedge v$ under the action of  $F_\Gamma^x$ (resp. $F_\Gamma$) generates $\Lambda^2T_xZ$ (resp. $\Lambda^2TZ$). 
\eprop

\bproof Recall that for every locally nilpotent vector field $\sigma$ on $X$ whose value at $x$ is $w \in T_xX$
and a function $f \in \Ker \sigma$ that vanishes at $x$ the phase flow of $f \sigma$ sends any vector $\nu\in T_xX$ into
the vector $\nu +t{\rm d}f(\nu )w$ where $t\in \C$ is the time parameter of the flow. 
Since $Z$ is $F$-flexible the fields $\sigma \in \cN$ generate $T_xZ$ and we can suppose that $u,v$, and $w$ are linearly independent (because  $\dim Z \geq 3$).
Thus choosing (by Lemma \ref{Gam}) a $\Gamma$-invariant $f$ so that ${\rm d}f(u )=0$ and $ {\rm d}f(v ) \ne 0$ we can transform $u\wedge v$ into $u \wedge v + u\wedge w$
under the action of the phase flow of a $\Gamma$-invariant replica of $\sigma$. That is,  we can get $u\wedge w$ in the span of the $F_\Gamma^x$-orbit of $u\wedge v$.
Applying the same argument one can get $w_1 \wedge w$ from $u \wedge w$.
Because  the fields $\sigma \in \cN$ generate $T_xZ$ we can view $w_1 \wedge w$ as a general element of $\Lambda^2T_xZ$ and thus
the orbit of $u\wedge v$ under the action of $F_\Gamma^x$ generates $\Lambda^2T_xZ$. 
Since $F$ acts transitively on $Z$ we have the desired conclusion.
\eproof

\section{Affine Homogeneous Spaces} \label{has}

\bnota\label{ahn1} From now on $G$ denotes a complex linear algebraic group and
$H$ is a proper closed algebraic subgroup of $G$. Let
$U$ (resp. $V$) be the unipotent radical of $G$ (resp. $H$).
That is, $G= M \ltimes U$ and $H= L \ltimes V$ where $M$ (resp. $L$) is a 
maximal reductive subgroup of $G$ (resp. $H$). By the Cartan-Iwasawa-Malcev theorem we can suppose further that $L$ is contained in $M$.
\enota


Recall the following (e.g., see \cite{Snow}).

\bthm\label{aht1} {\rm (1)} The homogeneous space $G/H$ is affine if and only if $G/V$ is affine.

{\rm (2)} Let $G$ be a unipotent linear algebraic group. Then $G/H\simeq \C^n$.
Moreover, if $G$ is solvable then $G/H \simeq \C^n \times (\C^*)^m$.
\ethm

\brem\label{ahr1}
Note that Theorem \ref{int1} implies AVDP for  $\C^n$ and $ \C^n \times (\C^*)^m$ and thus it suffices to consider the case when $G$ is not solvable. 
Similarly, we  suppose that $U$ is nontrivial since otherwise $G/H$ is affine if and only if $H$ is reductive (Matsushima's theorem), i.e. we have AVDP by Theorem \ref{int1}.
Furthermore, by analogous reasons we can suppose that $V$ is not contained in $U$ (see, \cite[Remark 6.3]{KaKu}). 
\erem








\bnota\label{ahn2}
Note that $M$ acts naturally on $G/V$ and since $M$ is reductive the quotient $Y=M\backslash G/V$ is affine when $G/V$ is such.
Furthermore, there is the following commutative diagram.

\bdi
G&\rTo^{V}&\, G/V \,\\
\dTo>{M}&&\dTo>{M}\\
U&\rTo^{V}&\, Y.\\
\edi

\vspace{1cm}

\enota

The next fact from \cite[Proposition 1]{Snow} elaborates it further.

\bprop\label{ahp1} If $G/H$ is affine then $G/V$ is $M$-equivariantly isomorphic to $M\times Y$ and $U$ is $V$-equivariantly
isomorphic to $Y \times V$.
\eprop 

\brem\label{ahr2} (1) Note that $U$ and $V$ are Euclidean spaces and therefore $Y$ is a factor in the cancellation problem. 



(2)  The variety $Y$  can be viewed as a submanifold of $U$ through identity $e \in U$ such that it is invariant under conjugation by any element of $L$.
More precisely, consider the natural $V$-equivariant morphism from $U \simeq Y \times V$ to $V$ (from Propostion \ref{ahp1}) and the conjugates $U \to V$ of this morphism
by elements of a maximal compact subgroup $L_\R$ of $L$.  Taking the average of these morphisms with respect to the Haar measure on $L_\R$
one gets an $L$-equivariant map $\kappa : U \to V$. Then $\kappa^{-1}(e)$ is the desired submanifold $Y$ (see the proof of Theorem 3 in \cite{Snow}).

In particular there is a natural action of $L$ on $Y$ induced by conjugations. Furthermore, the following is true \cite{Snow}. 
\erem

\blem\label{ahl0} The homogeneous space $G/H$ is isomorphic to $M\times_L Y$.
In particular there is a natural projection $\theta : G/H \to M/L$.
\elem

\bconv\label{ahconv1} (1) From now on we suppose that $X=G/H$ is affine.

(2) Let $S$ be a Levi subgroup of $M$, i.e. $M=S\ltimes T$ where $T$ is a torus. 
For every homogenous space $G/H$ one can suppose  without loss of generality (by going over to a covering), that $M=S\times T $, i.e., $S\cap T = \{ e \}$
which we shall do below. Furthermore, we suppose that $H$ and therefore $L$ are {\bf connected}.
\econv

\brem\label{ahr0} Of course, we need to consider more general affine homogeneous spaces $X'=G/H'$ where $H'$ is not necessarily connected but we can
always suppose that $H$ is a connected component of $H'$ while $L$ is a connected component of a maximal reductive subgroup $L'$ of $H'$
such that $L' \subset M$. Note that $L'=H'/V$ and $L=H/V$. Hence the finite group $\Gamma \simeq H'/H\simeq (H'/V)/(H/V) \simeq L'/L$ acts on $X=G/H$ and $X'=X/\Gamma$.
Furthermore, note that $\Gamma$ acts on $X$ from the right, i.e. this action commutes with the natural action of $G$ on $X$.

Recall also that in the presence of a $G$-invariant volume form $\omega$ on $G/H'$ to establish AVDP for $G/H'$ is the same as
to establish $\Gamma$-AVDP for $X$.
\erem

\bnota\label{ahn3}  Let  Notation \ref{ahn1} and \ref{ahn2} hold and  $S, X$ and $\Gamma$ be as in Convention \ref{ahconv1} and Remark \ref{ahr0},
i.e. $M=S \times T$ where $T$ is a torus.
Let 
$T_0$ be the image of $L$ under the natural homomorphism $M\to T$.
Then there exists a subtorus $T_1\subset T$ such that $T =T_0\times T_1$.  That is every $t \in T$ can be presented as $t=t_0t_1$ where $t_i \in T_i$
and by the Mostow theorem \cite{Mo} every $g\in G$ is of form $g=t_0t_1su$ where $s \in S$ and $u \in U$.
Consider the composition $\psi : G \to T$ of the natural quotient homomorphisms
$G \to M$ and $M\to T$ and the composition $\psi_1 : G \to T_1$ of $\psi$ and the natural projection $T \to T_1$. 
In particular $\psi (g)=t_0t_1$ and $\psi_1(g)=t_1$.

The composition of morphism $\theta : X \to M/L$ from Lemma \ref{ahl0} and the natural projection $M/L \to T_1$ will be denoted by
$\varphi : X \to T_1$.  

\enota

%

\blem\label{ahl1} The homomorphism $\psi_1$ coincides with the composition of $\rho : G \to G/H=X$ and $\varphi$.
Furthermore every fiber of $\varphi$ is an orbit of the natural action of the group $E=S \ltimes U$ on $X$, i.e. such a fiber $Z$
is a flexible variety (by \cite[Proposition 5.4]{AFKKZ}).
\elem

\bproof
Note that $\psi (L)=T_0$ and by the fundamental theorem of algebra $\psi (V)$ coincides with the identity $e \in T$.
Thus by the Mostow theorem $\psi (H)=\psi (LV)=T_0$.  Since $T_0 \subset \Ker \psi_1$ we see that the homomorphism $\psi_1$
(sending $g=t_0t_1su$ to $t_1$)
factors through $\varphi$ (sending the coset $gH=gL\cdot V=t_0t_1sL\tilde uV\in M \times_L Y$ to $t_1$ via $t_0t_1sL \in M/L$) which implies the first statement. 

Since $U$ is normal in $G$ and $T$ is in the center of $M$ the orbit $E g$ of $g$ can be written as $SUt_0t_1su=t_1t_0SU$.
Multiplying this orbit by $L$ from the right we get $EgL=t_1t_0SUL=t_1t_0SLU=t_1T_0SU$ since $SL=T_0S$. 
Note that $\rho$ sends each $H$-orbit and therefore each $L$-orbit to a point, i.e. $\rho (Eg)=\rho (EgL)$.
Since  $\psi_1^{-1}(t_1)=t_1T_0SU$ we see that $\rho (Eg)=\rho (\psi_1^{-1}(t_1))=\varphi^{-1}(t_1)$ which yields the second statement.
\eproof

\brem\label{ahr3} The first statement of Lemma \ref{ahl1} can be presented via the following commutative diagram. 

\bdi
G&\rTo^{}&G/U=M&\rTo^{}&T&\rTo^{}&T_1\\
\dTo>{}&&\dTo>{\simeq}&&\dTo>{\simeq}&& \dTo>{\simeq}\\
G/V=M\times Y&\rTo^{}&M&\rTo^{}&T&\rTo^{}&T_1\\
\dTo>{}&&\dTo>{}&&&& \dTo>{\simeq}\\
G/H=M\times_L Y&\rTo^{}&M/L&&\rTo^{}&&\, \, T_1 .\\
\edi

\vspace{1cm}
\erem

We finish this section  with some very special situation  which will be crucial for finding a locally nilpotent vector field  in Sections \ref{special case} and \ref{ADP}.

\bprop\label{ahp2} Let $M=S=SL_2$, $L \simeq \C^*$ be the group of diagonal matrices in $M$, the image of $H$ in $M$ under the natural morphism $\Theta : G \to M$
be the group of upper triangular matrices in $M$, and $V_0=V\cap U$ be different from $V$. Then
$Y \times \C$ is isomorphic to a Euclidean space. Moreover $Y$ admits a nontrivial locally nilpotent vector field $\eta$.

\eprop

\bproof 
Note that $\Theta (V)$ is a unipotent subgroup of the group of the  upper triangular $(2 \times 2)$-matrices, i.e. $\Theta (V) \simeq \C_+$. Since $\Ker \Theta |_V =V_0$
we have $V/V_0\simeq \C_+$. Now we can use the argument in the proof of \cite[Proposition 1]{Snow}. Namely, for  $G/V_0=M\times U/V_0$ one has the natural map $G/V_0 \to G/V$
which makes $G/V_0$ a bundle over $G/V$ with fiber $V/V_0$. Since $G/V$ is affine by Theorem \ref{aht1} this line bundle is trivial, i.e. $G/V_0 \simeq G/V \times \C$.
Consider the quotient $Y_0= M \setminus G /V_0 =U/V_0$. The triviality of first cohomology with coefficients
in the structure sheaf of the affine variety $Y$ implies that the principal $V/V_0$-bundle $Y_0 \to Y$ is again trivial, i.e. $Y\times \C \simeq Y_0 =U/V_0$
where the last variety is isomorphic to a Euclidean space by Theorem \ref{aht1}. This is the desired statement.
Since the ML-invariant of a Euclidean space is trivial the theorem contained in the papers of Makar-Limanov with Bandman and  Crachiola (\cite{BM-L}, \cite{CM-L})
implies that $Y$ admits a nontrivial locally nilpotent vector field $\eta$.

\eproof

\section{Condition (A)}

We keep Notation \ref{ahn3} in this section. 

\brem\label{mtIr1} Let $R=L\cap S$. By \cite[Lemma 4.7]{KaKu}
$M/L$ is $M$-equivariantly isomorphic to $S/R \times T_1$ and $M/L'$ is the quotient of $M/L$ with respect to the natural action of $\Gamma$.
Note that the multiplication by elements of any subgroup $F \subset S \subset M$ on $X=G/H=M\times_L Y$ from the left makes morphism
$X \to M/L$ in the last row of the diagram in Remark \ref{ahr3} $F$-equivariant and the action of $F$ on the base $M/L\simeq (S/R)\times T_1$ is nothing
but multiplication of cosets from $S/R$ by elements of $F$ from the left. Thus we have the following.

(1) If the action of $F$  on $S/R$ is nondegenerate (i.e. the dimension of a general orbit coincides with the dimension of $F$) and fixed point free then the natural
action of $F$ on $X$ is also nondegenerate and fixed point free.

(2) Since $F$ acts on the first factor of  $S/R \times T_1$ this action preserves
the fibers of $\varphi$, i.e. in (1) we have a nondegenerate and fixed point free action on each fiber $Z=\varphi^{-1}(t_1)$. 

(3) Multiplication of $X$ by elements of $T_1$ from the left generates a free action on the second factor of $M/L =S/R\times T_1$
which shows that morphisms $\theta : X \to M/L$ and $\varphi : X \to T_1$ are $T_1$-equivariant. In particular we have the following  extension of Lemma \ref{ahl1}.

\erem

\bprop\label{Anew} The homogeneous space $X$ is naturally isomorphic
to $Z\times T_1$ where $Z$ is a fiber of morphism  $\varphi$ (which therefore plays the role of the projection onto the second factor). 
\eprop

\bprop\label{ahp3} Let Convention \ref{ahconv1} hold and $X$ have a $G$-invariant volume form $\omega$.
Suppose that $X'=X/\Gamma$ is as in Remark \ref{ahr0}. Then $\chi_\omega$ is a character of $\Gamma$.

\eprop

\bproof By Remark \ref{ahr0} $\Gamma =L'/L$. Recall that  $\psi (L)=T_0$ and thus $\psi_1$ induces a well-defined homorphism $g : \Gamma \to \Gamma_1$ 
onto a finite subgroup $\Gamma_1$ of $T_1$. The action generated by left multiplication of points of $X$ by elements of $G$ commutes with the $\Gamma$-action.
Thus we can consider a new action of $\Gamma$ on $X$ which is obtained by composition of the original action of $\gamma \in \Gamma$ and
the left multiplication by $g(\gamma^{-1})$. Since the left multiplication
preserves $\omega$ the map $\chi_\omega$ is a character if and only if the new action of $\Gamma$ is also associated with a character.
Thus we switch to the new action whose advantage is in preservation of fibers of morphism $\varphi : X \to T_1$.
Left multiplications by elements of $T_1$ produce  commuting linearly independent semi-simple vector fields $\nu_1, \ldots , \nu_k$ (where $k=\dim T_1$) on $X$ 
which are $\Gamma$-invariant.
Taking interior products all of these fields with $\omega$ we get a volume  form $\omega (t_1)$ on a fiber $\varphi^{-1}(t_1)$. 
Since $\varphi^{-1}(t_1)$ is a homogeneous space of
the group $E$ from Lemma \ref{ahl1} it is rationally connected in terminology of \cite{KaKu4} (since the group $E$ is)
which implies that $\chi_{\omega (t_1)}$ is indeed a character by \cite[Proposition 6.4]{KaKu4}.

Hence any invertible function $\chi_{\omega} (\gamma )$ is constant on fibers of $\varphi$ and therefore it can be viewed
as a function on $T_1$. On the other hand since the $\Gamma$-action commutes with the left multiplications by elements of $T_1$
this function must be $T_1$-invariant and therefore constant. Thus we get the desired conclusion.

\eproof

\bnota\label{mtIn1}

Consider the subgroup $\SAut_{E}(X)$ of $\SAut (X)$ generated  by elements of  replicas of $\C_+$-actions associated with multiplications of elements of
$\C_+$-subgroups of $E=S \ltimes U$ (note that $\SAut_{E}(X)$ preserves each fiber $Z$ of $\varphi$ and  acts transitively on it, i.e. $Z$ is $\SAut_E(X)$-flexible). 
Let $\Gamma$ be the finite group acting on  $X$ from Remark \ref{ahr0}, i.e. every field associated with the natural left action of a $\C_+$ (or $\C^*$) subgroup
of $G$ is $\Gamma$-invariant.    Let $\SAut_E^\Gamma (X)$ be the subgroup of $\SAut_E (X)$ that commutes with
the action of $\Gamma$. Consider the set $\cW$ (resp. $\cW_1$) of wedge-products  $u\wedge v\in \Lambda^2T_zZ\subset \Lambda^2T_zX$
such that $u,v \in T_zX$ (resp. $u,v \in T_zZ$) be the values of module generating $\Gamma$-invariant vector fields on $X$ (resp. tangent to the fibers of $\varphi$).
\enota

By Proposition \ref{Gamma} 
we have the following.

\blem\label{mtIl1} Let $u\wedge v\in \cW_1$ be nonzero.
Then the span of the $\SAut_E^\Gamma (X)$-orbit of $u\wedge v$ contains $\Lambda^2T_zZ$.
\elem

Since at any point $t_1 \in T_1$ pairs of vectors that are values of semi-compatible vector fields
on the torus $T_1$ generate $\Lambda^2T_{t_1}T_1$ (indeed, consider the pairs $(\nu_i , \nu_j)$ where $\nu_i$'s are fields from the proof
of Lemma 4.2 in \cite{KaKu}) we have the following.

\blem\label{mtIl2}  Let $z \in Z$ and $t_1=\varphi (z)$. Then the
natural projection $X \to T_1$ induces a map  $ \Lambda^2T_zX\to  \Lambda^2T_{t_1} T_1$ such that the image of $\cW$ coincides with $ \Lambda^2T_{t_1} T_1$.

\elem

\bprop\label{A} Let for every fiber $Z$ of $\varphi$ the set $\cW_1$ from  Lemma \ref{mtIl1}  contains a  nonzero element.
Then Condition (A) of Theorem \ref{mct1} holds with module generating vector fields (see Remark \ref{mcr2}).
\eprop

\bproof
By Lemmas \ref{mtIl1} and \ref{mtIl2} it suffices to show that $\cW$ contains all wedge-products of form $u\wedge v$ where $u \in T_zZ$ and $v\in T_zX$ is such that
$u$ is a general vector of $T_zZ$ and
the image of $v$ in $T_{t_1}T_1$ (where $t_1=\varphi (z)$) is a vector of a basis of $T_{t_1}T_1$. Choose a $\C_+$ or $\C^*$-subgroup of $S$ associated
with a locally nilpotent vector field $\xi$ and a $\C^*$-subgroup in $T_1$ associated with a semi-simple vector field $\eta$. 
They are semi-compatible and $\Gamma$-invariant. Thus by Proposition \ref{scp1} they are $\Gamma$-semi-compatible,
and $\xi \wedge \eta$ yields an element $u \wedge v \in \cW$. 
Let us choose a general $\C_+$-subgroup in $E$ associated with a locally nilpotent vector field $\sigma$ whose value at $z$ is a general vector $w\in T_zZ$.
By Lemma \ref{Gam} we can choose a $\Gamma$-invariant replica $f\sigma$ such that its phase flow transforms $u\to u +tw$ and $v \to v +w_1$ in $T_zX$ where $t \ne 0$ and
the vector $w_1 \in T_zZ \subset T_zX$ because this phase flows preserves the fibers of $\varphi$.

Hence the span of the $\SAut_E^\Gamma (X)$-orbit of $u\wedge v$ contains $tw\wedge w_1 + u\wedge w_1+t w\wedge v$.
By Lemma \ref{mtIl1} $tw\wedge w_1 + u\wedge w_1 \in \Lambda^2 T_zZ$ is contained in $\cW$. Hence $w\wedge v \in \cW$ and we are done.
\eproof

\section{The validity of Condition (B)}

\bprop\label{B} Condition (B) of Theorem \ref{mct1}  holds for any affine $X=G/H$.

\eprop

\bproof Recall that $X=G/H=M\times_LY$ can be viewed as a locally trivial fibration with fiber $Y$ over $M/L=S/R\times T_1$
where $Y$ is contractible. In particular when $\dim Y \geq 2$ then the spectral sequence implies that $H^{n-1}(X)=0$ (where $n=\dim X$)
and Condition (B) holds automatically. If $\dim Y =0$ then $G$ is reductive and we are done by \cite[Proposition 4.5 and Lemma 4.11]{KaKu}.
Suppose that $\dim Y=1$, and thus $Y\simeq \C$. Then the desired conclusion follows from the next Proposition.

\eproof

\bprop\label{B1} Let $p : X \to B$ be a morphism of smooth affine algebraic varieties such that it is a locally trivial
$\AA^1$-fibration. 
Suppose that $X$ (resp. $B$) is equipped with a volume form $\omega$ (resp. $\omega_0$),
$\dim X=n$ (resp. $\dim B=n-1$),
and $\Theta : {\rm AVF} (X) \to \Omega^{n-1}(X)$ is the natural isomorphism between algebraic vector fields on $X$ and
algebraic $(n-1)$-forms induced by $\omega$. Then for every element of $H^{n-1}(X, \C)$ there exists a closed form $\tau \in \Omega^{n-1}(X)$ generating this element such
that it is the image
of a locally nilpotent vector field tangent to the fibers of $p$.
\eprop

\bproof  Since $X$ is a locally trivial $\AA^1$-fibration over $B$ we have a natural isomorphism $p^*: H^{n-1}(B, \C )\to H^{n-1}(X,\C)$.
In particular any nonzero element of $H^{n-1}(X, \C)$ can be presented by a closed form $\tau =p^*(\tau_0)$ where $\tau_0$ is a closed
$(n-1)$-form on $B$. Hence  for a general point $x \in X$ the form $\tau$ vanishes on a wedge-product $v_1 \wedge \ldots \wedge v_{n-1}$ of vectors $v_i \in T_xX$
if and only if a nonzero linear combination of the factors is tangent to the fiber of $p$ through $x$ (which is equivalent to the fact
that the vectors $p_*(v_1), \ldots , p_*(v_{n-1})$ are linearly dependent in $T_{p(x)}B$). Since the isomorphism $\Theta$ assigns to a vector field $\nu$ the interior product $\iota_\nu (\omega )$
we see that this property of $\tau$ implies that $\nu=\Theta^{-1}(\tau )$ is tangent to the fiber of $p$.
Recall that the restriction of $\Theta$ yields an isomorphism between the space ${\rm AVF}_\omega (X)$ of divergence-free (with respect to $\omega$)
algebraic vector fields and the space $\cZ_{n-1}(X)$ of closed $(n-1)$-forms. Note also that for a neighborhood $U$ of $b=p(x)\in B$ one has 
$p^{-1}(U)\simeq U \times \C$ (where the second factor has a coordinate $z$) and
$$\omega |_{p^{-1}(U)}=g\omega_0|_U \wedge {\rm d} \, z$$ where $g$ is an invertible function on $U$.
This implies that the field $\nu =f \partial /\partial z$ (where $f$ is a function on $p^{-1}(U)$) has a zero $\omega$ divergence if and only
if the restriction of $\nu$ to $p^{-1}(b)\simeq \C_z$ has the zero divergence with respect to 1-form ${\rm d} \, z$. The only vector fields
on the line $\C$ with this property are locally nilpotent vector fields  $a \partial / \partial z$ with $a \in \C$. Thus $\nu |_{p^{-1}(U)}$ is of form
$h \partial / \partial z$ where $h$ is a function on $U$ and $\nu$ is locally nilpotent which is the desired conclusion.

\eproof

\section{AVDP for affine homogeneous spaces, the general case}

In this section we consider affine homogeneous space $X=G/H$ such that under Convention \ref{ahconv1}  the quotient $S/R$ is not isomorphic to $SL_2/\C^*$ where $\C^*$
is the subgroup of diagonal matrices in $SL_2$.

The reason why this case is rather simple is due to the following fact \cite[Theorem 24]{DDK}.

\bprop\label{DDK} Let $S$ be non-trivial and $S/R$  be non-isomorphic to ${SL}_2/\C^*$.
Then there is a pair locally nilpotent compatible (and in particular semi-compatible) vector fields $\xi$ and $\eta$ on $X$ tangent to each fiber $Z=\varphi^{-1}(t_1)$
whose phase flows are associated with $\C_+$-subgroups of some $F \simeq SL_2$ in $S$. 
\eprop

Hence applying  Proposition \ref{A} we have the following.

\bcor\label{newcor} Let $S$ be non-trivial and $S/R$  be non-isomorphic to ${SL}_2/\C^*$. Then condition (A) from Theorem \ref{mct1} is valid.

\ecor

Corollary \ref{newcor} and Proposition \ref{B} imply now that all assumptions of  
Theorem \ref{mct1} are valid in the special case of nontrivial  $S$ and $S/R$ non-isomorphic to $SL_2/\C^*$. Since the case of a  trivial $S$ was dealt in Remark \ref{ahr1}
we have the following part of our main result.

\bcor\label{dim3} Let an affine  $X=G/H$ have a left-invariant volume form $\omega$ and assume $S/R$ is not isomorphic to $SL_2/\C^*$. 
Then $X$ has AVDP with respect to $\omega$.
\ecor

\section{$\C [X]$-modules in $\AVF (X)$}
\label{newtech1}

The case of $S/R \simeq SL_2/\C^*$ is more complicated not only for establishing AVDP but also for ADP. The next two facts will help
us to deal with the latter.

\bthm\label{ap1} Let $\pi : X \to B$ be a surjective morphism of smooth affine algebraic varieties, i.e. $\C [B]$ is naturally embedded
into $\C [X]$. Suppose that $\eta$ is a locally nilpotent or  semi-simple vector field on $X$ such that $\C [B] \subset \Ker \eta$.
Let $\alpha_i, \, i=1, \ldots , k$ be regular functions on $B$ without common zeros and $D_i=( \alpha_i\circ \pi)^{-1}(0)$ (in particular $\bigcap_{i=1}^kD_i = \emptyset$).
Suppose that $\xi_i, \,  i=1, \ldots , k $ are locally nilpotent vector fields on $X$ such that 
$\alpha_i \in \Ker \xi_i$ (and thus $\xi_i$ is tangent to $D_i$) and
the restriction of the pair $(\xi_i, \eta)$ is compatible on $X_i=X \setminus D_i$.

Then there exists a nontrivial ideal $J\subset \C [X]$ for which $J \eta \subset \LieA (X)$.
\ethm

\bproof Let $I_i \subset \C [X_i]$ be the associated ideal of the compatible pair $(\xi_i , \eta )$. That is, every function $h \in I_i$ can be 
presented in the form $h=\sum_{j=1}^n f_jg_j$ with $f_j \in \Ker \xi_i\subset \C [X_i]$ and $g_j \in \Ker \eta \subset \C [X_i]$.
Note that for every $f_j$ there exists a natural $m_j = m_j (f_j)$ for which $\alpha_i^{m_j} f_j \in \C [X]$ (and similarly $\alpha_i^{n_j} g_j \in \C [X]$ for some $n_j$).
In particular for every $h \in I_j$ there exists a natural $l_i=l_i(h)$ so that $\alpha_i^{l_i} h = \sum_{j=1}^n f_jg_j$ with $f_j \in \Ker \xi_i\subset \C [X]$ and $g_j \in \Ker \eta \subset \C [X]$.
Let  $b_i =\xi_i (a_i)$ where $a_i \in \C [B]\setminus \Ker \xi_i$ plays the role of function $a$ from Definition \ref{mcd1} (that is,
$\xi_i^2 (a_i)=0$). Using the formula
$$ [f_j\xi_i, a_ig_j\eta ] - [a_i f_j\xi_i, g_j\eta] = b_if_jg_j \eta$$
we see that that the field $\alpha_i^{l_i} b_i h \eta $ is contained in $\LieA (X)$. Furthermore, replacing $g_j$ by $\beta_i g_j$ we see that
 $\alpha_i^{l_i} \beta_i b_i h \eta $ is contained in $\LieA (X)$ for every $\beta_i \in \C [B]$. 
 Thus for a given $i$ and every element $e=b_ih$ of the ideal $L_i=b_iI_i$  the field $\alpha_i^{l_i} \beta_i e \eta $ is contained in $\LieA (X)$. 
 
Denote by $J_i$ the ideal $\C [X] \cap L_i$ and let $J = \bigcap_{i=1}^k J_i$. That is, for every $e\in J$ the field 
$\alpha_i^{l_i} \beta_i e \eta $ is contained in $\LieA (X)$ already for every index $i=1, \ldots , k$.
Choosing (by the virtue of the Nullstellensatz) $\beta_j$ so that $\sum_{j=1}^k\alpha_i^{l_i}\beta_i=1$ we see that
$e \eta $ is contained in $\LieA (X)$ which concludes the proof.
\eproof

\bcor\label{ap2}  Let  $X$ and $B$ be smooth affine algebraic varieties and let
$\pi : X \to B$ be a Zariski locally trivial fibration with fiber $Y$ admitting a nontrivial locally nilpotent vector field $\delta$.
Let $\alpha_i, \, i=1, \ldots k$ be regular functions on $B$ without common zeros such that $X_i:=X\setminus \alpha_i \circ \pi^{-1}(0)$ is naturally isomorphic to
the direct product $X_i=B_i\times Y$ where $B_i=B\setminus \alpha_i^{-1}(0)$.
Suppose that $\xi_i, \, i=1, \ldots , k$ are locally nilpotent vector fields on $X$ such that 

{\rm (i)} $\xi_i (\C [B]) \subset \C [B]$, i.e. $\xi_i$ induces a locally nilpotent vector field $\xi_i'$ on $B$;

{\rm (ii)} $\alpha_i \in \Ker \xi_i'$, i.e. $\xi_i'$ is locally nilpotent on $B_i$;

{\rm (iii)} $\xi_i'$ admits a slice\footnote{That is, $B_i$ is isomorphic to $R_i \times \C$ and the action
is nothing but a translation on the second factor.}, and in particular it is nontrivial.

Then there exist a nontrivial locally nilpotent vector field $\eta$ and a nontrivial ideal $J\subset \C [X]$ for which $J \eta \subset \LieA (X)$.
\ecor

\bproof Note that $\delta$ induces a locally nilpotent vector field on $X_1$ whose kernel contains $\C [B\setminus \alpha_1^{-1}(0)]$.
Multiplying this field by some power of $\alpha_1$ we get a locally nilpotent vector field $\eta$ on $X$ with $\C [B] \subset \Ker \eta$.
Condition (i) implies that the $\C_+$-action induced by $\xi_i$ (resp. $\xi_i')$) on $X$ (resp. $B$) makes $\pi$ equivariant. Furthermore the existence of
a section $R_i\subset B_i$ for $\xi_i'$ implies that $\pi^{-1}(R_i) $ is a section of $\xi_i|_{X_i}$, and thus 
$X_i$ is isomorphic to $\pi^{-1}(R_i) \times \C$ with $\xi_i|_{X_i}$ being a translation on the second factor.
That is, the image of $\C [\pi^{-1} (R_i)]$ in $\C [X_i]$ under the natural embedding is contained in the kernel of $\xi_i|_{X_i}$. Since $\pi^{-1}(R_i) \simeq R_i \times Y$
we see that a trivialization $X_i \simeq B_i\times Y $ can be chosen so that $\xi_i|_{X_i}$ has $\C [Y]$ in its kernel. 
Hence $\xi_i$ is semi-compatible with $\eta$. Furthermore, since $\xi_i$ is nontrivial on $B$ there exists an
element $a_i \in \C [B]\setminus \Ker \xi_i$ for which $\xi_i^2 (a_i)=0$. That is, the pair $(\xi_i, \eta)$ is compatible on $X_i$
and the desired conclusion follows from Theorem \ref{ap1}.

\eproof

\section{$\C [X]$-modules in $\cC_{n-2}(X)$}
\label{newtech2}

To deal with the difficult special case of $S/R \simeq SL_2/\C^*$ we need the next two facts which guarantee the existence of module generating pairs different from  semi-compatible pairs.

\bnota\label{ap3}  Suppose that a smooth affine algebraic variety $X$ of dimension $n$ is equipped with an algebraic volume form $\omega$.
Let    $ \AVF_\omega (X)$, $\cC_{k} (X)$,    $\cZ_{k}(X)$ and $\cB_{k} (X)$ be  as in Notation \ref{mcn2}. Recall that in this case 
there exists an isomorphism $\Theta : \AVF_\omega (X) \to \cZ_{n-1}(X)$
given by the formula $\xi \to  \iota_\xi \omega$ where $ \iota_\xi \omega$ is the interior product of $\omega$ and $\xi \in \AVF_\omega (X)$.

Furthermore for any two fields $\sigma$ and $\delta$ from $\AVF_\omega (X) $  one has the following \cite{KaKu}
\be\label{ape} \iota_{[\sigma, \delta]}={\rm d} \, \iota_\sigma \iota_\delta \omega \ee
where ${\rm d} : \cC_{n-2}(X) \to \cB_{n-1}(X)$ is the differential map.

\enota

\bthm\label{ap4} Let Notation \ref{ap3} hold and $\pi : X \to B$ be a surjective morphism of smooth affine algebraic varieties, i.e. $\C [B]$ is naturally embedded
into $\C [X]$. Suppose that $\eta \in  \IVF_\omega (X) $ is such that $\C [B] \subset \Ker \eta$.
Let $\alpha_i, \, i=1, \ldots , k$ be regular functions on $B$ without common zeros and $D_i=( \alpha_i\circ \pi)^{-1}(0)$.
Suppose that $\xi_i, \,  i=1, \ldots , k $ are 
fields from $\IVF_\omega (X)$ such that 
$\alpha_i \in \Ker \xi_i$ (and thus $\xi_i$ is tangent to $D_i$) and the pair $(\xi_i , \eta)$ is semi-compatible on $X_i = X \setminus D_i$.
Suppose also that

{\rm ($*$)} these fields $\{ \xi_i \}_{i=1}^k$ generate a $\C [X]$-module $N$ such that 
for every point $x \in X$ the localization of $N$ at $x$ coincides with the localization of the module generated by the fields $\{ \xi_j | \alpha_j (x) \ne 0 \}$. 

Then ${\rm d}^{-1} (\Theta (\LieAO (X))$ contains a $\C [X]$-submodule of $\cC_{n-2} (X)$ that consists of all $(n-2)$-forms that can be presented
as $$\tau = \iota_\zeta \iota_\eta \omega$$ with $\zeta$ running over $JN$ for some nontrivial ideal $J$ in $\C [X]$.

\ethm

\bproof Let $I_i \subset \C [X_i]$ be the associate ideal of the semi-compatible pair $(\xi_i, \eta )$. That is, every $h \in I_i$ can
be presented as $h=\sum_j f_jg_j$ where $f_j \in \Ker \xi_i$ and $g_j\in \Ker \eta$.  As in the proof of Theorem \ref{ap1} we note
that for some natural $l_i=l_i(h)$ one has $\alpha_i^{l_i} h = \sum_{j=1}^n f_jg_j$ with $f_j \in \Ker \xi_i\subset \C [X]$ and $g_j \in \Ker \eta \subset \C [X]$.
Note that $f_j\xi_i$ and $g_j\eta$ are still fields from  $\IVF_\omega (X)$.

Recall that for a differential form $\tau$, a regular function $f$, and a vector field $\zeta$  on $X$ one has the following property of the interior product: 
$\iota_{f\zeta}\tau =f \iota_\zeta \tau$. Hence
letting $\sigma =f_j\xi_i$ and $\delta =g_j\eta$ in Formula (\ref{ape}) we see that $f_jg_j  \iota_{\xi_i} \iota_\eta \omega$ is contained
in ${\rm d}^{-1} (\Theta (\LieAO (X))$, and therefore
\be\label{ape2} \alpha_i^{l_i} h  \iota_{\xi_i} \iota_\eta \omega \in {\rm d}^{-1} (\Theta (\LieAO (X)). \ee
In the case when $h$ belongs to the ideal $J=\bigcap_{i=1}^k (I_i \cap \C [X])$ Formula (\ref{ape2}) holds for every index $i$.

Consider the set $N'$ of all linear combinations of forms from Formula (\ref{ape2}) with $h$ running over $J$ and $i$ changing from 1 to $k$.
Replacing every $g_j$ by $\beta g_j$ with $\beta \in \C [B]$  (i.e. $\beta g_j$ is still in $\Ker \eta$) we see that for every $\tau \in N'$ one has $\beta \tau \in N'$; that is, $N'$ is a $\C [B]$-module.
Furthermore, every element of $N'$ can be presented as $\iota_\zeta \iota_\eta$ where $\zeta$ is a vector field of 
form $\zeta = \sum_{i=1}^k  \alpha_i^{l_i} h_i  {\xi_i} $ with $h_i \in J$. Such fields $\zeta$ form a $\C [B]$-module $L$. 

We need to prove the equality $L=JN$. Let $\cN$ and $\cL$ be quasi-coherent sheaves over $B$ generated by $N$ and $L$ respectively. 
Since by \cite[Chapter 2, Corollary 5.5]{Har} there is an equivalence between categories of $\C [B]$-modules and of quasi-coherent sheaves over $B$
it suffices to check the equality $\cL = J \cN$. Note that if $\alpha_i (x) \ne 0$ then  the fact $  \alpha_i^{l_i} h  {\xi_i} \in L$ implies that
$h\xi_i \in \cL$. Furthermore, by  Condition ($*$)  $J \cN$ is generated by the germs of such vector fields. 
This implies the equality $\cL = J \cN$ and hence the desired conclusion.
\eproof

\bcor\label{ap5}  Let  
$\pi : X \to B$, $Y$, $\delta$, $\alpha_i, \, i=1, \ldots k$,  $B_i$, and $X_i:=X\setminus \alpha_i \circ \pi^{-1}(0) \simeq B_i\times Y$
be as in Corollary \ref{ap2}.  
Suppose that $\pi$ is equivariant with respect to some actions of a semi-simple group $S$ on $X$ and on $B$. Let  $\xi_i, \, i=1, \ldots , k$ be
a collection of locally nilpotent vector fields on $X$ induced the actions of 
one-parameter unipotent subgroups of $S$  such that  conditions (i)-(iii) from Corollary \ref{ap2} are valid.  Let $N$ be the $\C [X]$-module generated by $\{ \xi_i \}_{i=1}^k$.
Suppose that

{\rm ($\star$)} for every $x \in X$ the subspace $N_x$ of $T_xX$ generated by the values of $N$ at $x$
coincides with the subspace generated by the values of the fields $\{ \xi_j | \alpha_j (x) \ne 0 \}$. 

Then in the case when $X$ is equipped with a volume form $\omega$ the set ${\rm d}^{-1} (\Theta (\LieAO (X))$ contains a $\C [X]$-submodule of $\cC_{n-2} (X)$ that consists of all $(n-2)$-forms that can be presented
as $$\tau = \iota_\zeta \iota_\eta \omega$$ with $\zeta$ running over $JN$ for some nontrivial ideal $J$ in $\C [X]$.

\ecor

\bproof As in the proof of Corollary \ref{ap2} we see that $\delta$ generates  a locally nilpotent vector field $\eta$ on $X$ with $\C [B] \subset \Ker \eta$.
Hence the  semi-compatibility of each $\xi_i$
with $\eta$ follows  as before from the fact that a trivialization  $X_i \simeq B_i\times Y $ can be chosen
so that the restriction of $\xi_i$ to $X_i$ has $\C [Y]$ in its kernel. 
Note also that Condition ($\star$) implies Condition ($*$) in Theorem \ref{ap4} which yields the desired conclusion.

\eproof

\section{AVDP for affine homogeneous spaces, the special case}\label{special case}

\bnota\label{mtIIn1} We keep Notation \ref{ahn3} in this section.
Until Theorem \ref{main} we suppose that $M=S\times T$ where
$S=SL_2$ and $R=L\cap S$ is a proper connected nontrivial reductive subgroup in $SL_2$, i.e.  $R \simeq \C^*$.
 Furthermore we suppose that $V$ is not contained in $U$ (because otherwise by Remark \ref{ahr1} we are in a case considered before).
\enota

\blem\label{mtIIl1} Let $H'$ and $\Gamma$ be as in Remark \ref{ahr0}. Then $\Gamma$ is trivial, i.e. $H'=H$.

\elem

\bproof 
Recall that for under the natural homomorphism $\psi : G \to T$ one has $\psi (L)=T_0$ and $\Ker \psi |_L =R$.
Hence $L=R \times T_0$ (since $M=S \times T$).  Suppose that $L'$ is as in Remark \ref{ahr0} and $R'=L' \cap S$. 
Consider an element $l'=r't' \in L' $ where $r' \in S$ and $t' \in T$.
Multiplying $l'$ by an element of $T_0$ we can suppose that $t' \in T_1$.
Since $L$ is a normal subgroup of $L'$ (i.e. $(l')^{-1}Ll'=L$) we see that $r'$ belongs to the normalizer $N_S(R)$ of $R$ in $S\simeq SL_2$ and since $R$ can be viewed
as the subgroup of diagonal matrices in $SL_2$ we see that
$N(R)/R \simeq \Z_2$. Hence $l'$ is contained
in $N_S(R)T_1$. The subgroup $\psi (L')\cap T_1$ is finite and replacing $T_1$ by the torus $T_1/(\psi (L')\cap T_1)$
we can suppose that it is trivial. Hence $l'\in N_S(R)$ and simultaneously $l' \in R'=L'\cap S$. Thus $\Gamma = L'/L$
is isomorphic to $R'/R$ and depending on whether $R'=R$ or $R'=N_S(R)$ we get $\Gamma$ either trivial or isomorphic to $\Z_2$.
Thus it suffices to show that $R'=R$.

Let $V_0=V \cap U$ and let
$L_0$ be the image of $V$
under the natural homomorphism $G \to G/U=M=S\times T$. By Notation \ref{mtIIn1} $V_0 \ne V$ and $L_0$ is nontrivial, and it is also unipotent. 
In particular it must be contained in $S$ and since the only  nontrivial unipotent subgroup of $SL_2$ is $\C_+$ we suppose that
$L_0\simeq \C_+$. 





Since $V_0$ is a normal subgroup of $H$ one has $rV_0r^{-1}=V_0$  for every $r \in R$ (resp. $r \in R'$) .  Applying the natural
morphism $G \to S$ we get $rL_0r^{-1}=L_0$. 

Therefore the Lie algebra of $L_0$
consists of eigenvectors for the adjoint action of $R$. Treating $R$ as the group of diagonal matrices in $SL_2$ we see that
the only eigenvalues of the adjoint action on $\cS\cL_2$ are matrices with only nonzero entries in the upper right (resp. lower left) 
corner.  
Hence
we can also treat $L_0$ as the subgroup of the upper unipotent matrices and $R$ as the subgroup of the diagonal matrices.
However for $r' \in N_S (R) \setminus R$ the equality
$rL_0r^{-1}=L_0$ does not hold, i.e. $R=R'$. Thus we have the desired conclusion.
\eproof

We need the following technical fact.

\blem\label{mtIIl2}
Let $R\simeq \C^*$ (resp. $F\simeq \C_+$) be the group of diagonal (resp. lower unipotent) matrices $$\lambda= \left[
\begin{array}{rr}
a & 0  \\
0 & a^{-1} \\
\end{array}  \right] \, \, \,  (resp. \, \, \, \alpha =\left[
\begin{array}{rr}
1 & 0  \\
t & 1 \\
\end{array}  \right])$$
in $SL_2$.
Then

 {\rm (1)} the homogeneous space $K=SL_2/R$ is naturally isomorphic to the surface $ \{ xy=z^2 -z\} \subset \C^3_{x,y,z}$;

{\rm (2)} the left multiplication by elements $\alpha$ of $F$
induces the $\C_+$-action on $K$ given by $$(x,y,z) \to (x, y+(2z-1)t +xt^2, z+xt).$$  
In particular, this action preserves subvarieties $K^i=K \setminus \{ x=z-i=0 \}\simeq \C_{x,y}^2$ for $i=0,1$. 

\elem

\bproof Consider the element
\begin{center}
$ \left[
\begin{array}{rr}
a_1& a_2  \\
b_1 & b_2 \\
\end{array}  \right] \in SL_2$ \end{center}
and note that the right multiplication by $\lambda$ yields the matrix
\begin{center}
$ \left[
\begin{array}{rr}
a a_1& a^{-1}a_2  \\
a b_1 & a^{-1} b_2 \\
\end{array}  \right] .$ \end{center}
Note that functions $x:=a_1a_2, y:=b_1b_2$, and $z:=a_1b_2$ generate the invariant subring of $\C [Y]$ (with respect the 
the $\C^*$-action on $K$ associated with right multiplications by elements of $R$) and the equality $a_1b_2-a_2b_1=1$ yields
$xy=z^2-z$. The rest is the direct computation.
\eproof

\bnota\label{mtIIn2} Consider functions $x,y,z$ on $K$ from Lemma \ref{mtIIl2}. They generate functions on $M/L=K\times T_1$
and on any other variety over $K$ which
by abuse of notation we denote by the same symbols. 
\enota

\blem\label{mtIIl4} Let $\tilde K^i=\theta^{-1} (K^i\times T_1)$ (recall that by Remark \ref{mtIr1} $M/L \simeq S/R \times T_1=
K\times T_1$ where $K$ is  from Lemma \ref{mtIIl2}).
Let $\cS^0$ (resp. $\cS^1$) be given by equations $y=z-1=0$ (resp. $y=z=0$)
in $M/L$, i.e. $\tilde \cS^i := \theta^{-1}(\cS_i)$ is contained in  $\tilde K^i$.
Then

{\rm (1)}  there are natural isomorphisms $\tilde \cS^i\simeq \cS^i \times Y$ and $\tilde K^i \simeq \C_{y} \times  \tilde \cS^i \simeq \C^2\times Y$,
such that the restriction of   the action of the group $F\subset SL_2$ (of lower unipotent matrices) to $\tilde K^i$ is a translation along the factor $\C_y$ (i.e. $\tilde \cS^i$ is a slice of this action);

{\rm (2)} morphism $\theta : X \to M/L$ is a Zariski locally trivial fibration with fiber $Y$;

{\rm (3)} For every point $p \in X$ there an element of $s \in SL_2$ such that after the left multiplication by $s$ one
has $s.p \in \tilde K^0$.

\elem

\bproof 

Exchanging the role of $x$ and $y$ in Lemma \ref{mtIIl2} one gets $\C_+$-action on $K$ given by $$(x,y,z) \to (x+(2z-1)t +yt^2, y, z+xt).$$
Since it is induced by left multiplications by the group $F'$ of upper unipotent matrices it generates also generates a $\C_+$-action on $X$ such that $\tilde \cS^i$ is invariant under it.
For isomorphism  $\tilde \cS^i\simeq \cS^i \times Y$ it remains to note that, say, 
$\{ x=y=z-1=0\} \subset \tilde \cS^0$ is isomorphic to $Y$ and it is a slice of the restriction of the $\C_+$-action on $\tilde \cS^0$. 

For isomorphism $\tilde K^i \simeq \C_{y}\times \tilde \cS^i$ it suffices to note that
$\tilde \cS^i$ is a slice of the $F$-action on $\tilde K^i$ because, say, the line $\{ y=z-1=0\} \subset K$ is a slice of the $F$-action on the surface $K^0\simeq \C^2$.
Thus we have (1).

Statement (2) follows from the fact that $K^0 \cup K^1$ coincides with $K$.

Furthermore, $K$ is transitive with respect to the group of algebraic automorphisms generated by elements of the $F$ and $F'$-actions (e.g., see \cite{ML90}). This yields (3).

\eproof

Now we can prove our main result.

\bthm\label{main} Let $X$ be an affine homogeneous space of a linear algebraic group $G$. Suppose that $X$
has a left-invariant (with respect to the natural $G$-action) volume form $\omega$.
Then $X$ has AVDP with respect to $\omega$.

\ethm

\bproof By Corollary \ref{dim3} it suffices to consider the case when $S=SL_2$ and $R\simeq \C^*$.
Furthermore by
Lemma \ref{mtIIl1} $\Gamma$ is trivial and $H=H'$. 

Thus in order to establish AVDP we need to check Conditions 
(A) and (B) of Theorem \ref{mct1} with trivial $\Gamma$-action on $X=G/H$. Condition (B) is valid by Proposition \ref{B} and
to prove (A) one has to find a nonzero element in $\cW_1$ from Lemma \ref{mtIl1}, i.e. to establish existence of a module generating
pair of vector fields on $Z$ where $Z$ is a fiber of $\varphi : X \to T_1$. 

By Lemma \ref{mtIIl4} $\theta : X \to M/L$ is a locally trivial Zariski fibration with fiber $Y$ is required in Corollary \ref{ap2}.
Consider the action of the group $F$ from Lemma \ref{mtIIl4} and a locally nilpotent vector field $\xi$ on $X$
generated  by it. Since $\theta$ is equivariant with respect to this action the field $\xi$ satisfies Condition (i) from Corollary \ref{ap2}
while Conditions (ii) and (iii) follow from Lemma \ref{mtIIl4} (with function $x$ playing the role of $\alpha_i$).

Conjugating the $F$ action by an element $s \in SL_2$ we get another locally nilpotent vector field $\xi_j$ on $X$ with similar properties
but such that by Lemma \ref{mtIIl4} the corresponding function $\alpha_j$ does not vanish at a given point $p \in X$.
Hence we can choose a finite number of such fields with no common zeros for functions $\alpha_j$. Hence all assumptions of Corollary \ref{ap2} are true.

Furthermore, conjugation by elements of $SL_2$ enables us to choose
this collection $\{ \xi_j \}$ of fields so that at every point $p \in X$ those of them that have $\alpha_j (p) \ne 0$
generate the same subspace of $T_pX$ as the whole collection. That is, Condition ($\star$) from Corollary \ref{ap5} is satisfied.
Hence for $\xi$ (as for any $\xi_j$) the pair $(\xi ,\eta )$ is module generating where $\eta$ is as in Corollary \ref{ap5}. Hence $\cW_1$ is not empty and we are done.


\eproof

\section{ADP for affine homogeneous spaces}\label{ADP}

The following fact can be found in \cite{KaKu1}.

\bprop\label{adpp1} Let $\{ (\xi_j,\eta_j) \}$ be a collection of compatible pairs of vector fields on $X$
and let $I_j$ be the associated ideal of the pair $ (\xi_j,\eta_j) $. Suppose that
that for every $x \in X$ the values of $\{ I_j(x)\eta_j \}$ at this point generate $T_xX$ where $I_j (x)$ is the set of values of $I_j$ at $x$.
Then $ \Lie (X) =\AVF (X)$, i.e. $X$ has ADP.
\eprop

\brem\label{adpr1}
(1) Note that if $X$ is homogeneous space that it suffices to require that the values  of $\{ \eta_j \}$ generate $T_xX$
only for one point $x \in X$. Furthermore, if the isotropy subgroup of $x$ in the group of algebraic automorphisms 
acts irreducibly on $T_xX$ it is enough to find one compatible pair $(\xi , \eta )$ such that the field $I\eta $ does not vanish at $x$.

(2) It is worth mentioning that similarly to Remark \ref{mcr2} the requirement about compatibility of pairs  $(\xi_j,\eta_j)$ is needed only to guarantee
that  $\LieA (X)$ contains the $\C [X]$-modules $I_j\eta$.
Thus one can replace the compatibility requirement in
Proposition \ref{adpp1} by the assumption that $I_j \eta$ is contained in $\LieA (X)$.
\erem

\bnota\label{adpn1} The set of complete $\Gamma$-invariant algebraic vector fields will
be denoted by $\IVF (X, \Gamma)$ and the Lie algebra generated by it by $\LieA (X, \Gamma)$. 
\enota

\bdefi\label{adpd1}
Let $X$ be a $\Gamma$-variety. Note that the fact that $X/\Gamma$ has ADP is equivalent to the fact
that $\LieA (X, \Gamma ) = \AVF (X, \Gamma )$. In the case of the last equality we say that $X$ has $\Gamma$-ADP.

\edefi

Thus we have the following.

\bprop\label{adpp2} {\rm (1)} Let $X$ be a $\Gamma$-variety homogeneous with respect to an action of linear algebraic group $G$ which
commutes with the action of $\Gamma$. Suppose that $\{ (\xi_j,\eta_j) \}$ is a collection of $\Gamma$-compatible pairs of vector fields on $X$
such that at some $x \in X$ the values of $\{ \eta_j \}$ at this point generate $T_xX$. Then $X$ has $\Gamma$-ADP.

{\rm (2)} Furthermore, if in (1)  the isotropy subgroup of $x$ in the group of algebraic automorphisms that
commute with the $\Gamma$-action acts irreducibly on $T_xX$ then for $\Gamma$-ADP of $X$  it suffices to find 
a $\Gamma$-compatible pair $(\xi , \eta )$ such that the field $\eta$ does not vanish at $x$.

\eprop

\blem\label{adpl1} Let $W=X\times Y$ where $X$ is as in Proposition \ref{adpp2} (2) and $Y$ and $W$ be $\Gamma$-varieties for which  the projections
$W \to X$ and $W \to Y$ are $\Gamma$-equivariant. Suppose that  $(\xi , \eta )$ is as in Proposition \ref{adpp2} (2)
(i.e. $X$ has $\Gamma$-ADP) and at any point $y\in Y$ the values of
$\Gamma$-invariant complete algebraic vector fields generate $T_y Y$. Then $W$ has $\Gamma$-ADP. 
\elem

\bproof Consider a point $w:=(x,y)\in W$ and the tangent space $T_wW = T_xX \oplus T_yY$. By Proposition \ref{adpp2} we have
$\Gamma$-compatible pairs $\{ (\xi_i , \eta_i )$ on $X$ (generated by application of the isotropy group to $(\xi , \eta )$)
such that the fields $\{ \eta_i \}$ generate
$T_xX \subset T_wW$. Consider now any $\Gamma$-invariant complete field $\nu$ on $Y$ and
the induced pair $(\xi , \nu )$ on $W=X\times Y$. By Example \ref{mce1} it is $\Gamma$-compatible which implies the desired conclusion.

\eproof

\bthm\label{main1}
Let  $X'$ be  a connected affine homogeneous space of a linear algebraic group $G$.
Suppose that $X'$ is not isomorphic to a line or a torus $(\C^*)^k$. Then $X'$ has ADP.
\ethm

\bproof As in Remark \ref{ahr0} we can suppose that $X' =X /\Gamma$ where $X$ satisfies Convention \ref{ahconv1}.
Recall that by Proposition \ref{Anew} (3) $X$ is isomorphic to $Z \times T_1$ and the $G$-action on $X$ commutes with the $\Gamma$-action.
Note also that the $\Gamma$-action preserves the fibers of the natural projections $X \to Z$ and $X \to T_1$, i.e. we can view
these projections as $\Gamma$-equivariant. 
Thus by Proposition \ref{adpp2} and Lemma \ref{adpl1} it suffices to find a $\Gamma$-compatible pair $(\xi , \eta )$ on 
the fiber $Z$ of morphism $\varphi : X \to T_1$ (recall that $Z$ is a homogeneous
space of the group $E =S\ltimes U$). As in Propositions \ref{Gam} and \ref{Gamma} we can check that the isotropy subgroup
of a given point $z\in Z$ in the group of automorphism that commute with $\Gamma$ acts irreducibly on $T_zZ$.
Hence if $S/R$ is not isomorphic to $SL_2/\C^*$ (as in Proposition \ref{DDK}) the existence of such a $\Gamma$-compatible pair
follows from Proposition \ref{DDK}. 

Thus consider the case when $M=S=SL_2$ and $L\simeq \C^*$ (recall that then $\Gamma$ is
trivial by Lemma \ref{mtIIl1}). As we mentioned  $X$ is naturally isomorphic to $Z\times T_1$.
Since at every point $t \in T_1$ complete algebraic vector fields generate the tangent space of $T_1$ at $t$, by Lemma \ref{adpl1} it suffices to check
ADP for $Z$ which is a Zariski locally trivial $Y$-fibration over $S/R$. Recall that $Y$ admits a nontrivial locally nilpotent vector field
by Proposition \ref{ahp2}.  Thus by Corollary \ref{ap2} there is a nontrivial locally nilpotent vector field $\eta$ on $Z$ for which
$\C [S/R]$ is contained in $\Ker \eta$. 

Furthermore, as we checked in the proof of Theorem \ref{main} there is a collection $\{ \xi_j \}$ of locally nilpotent vector fields 
and functions $\{\alpha_j \}$ satisfying all assumptions of Corollary \ref{ap2}. Hence $\LieA (X)$ contains all field of form
$J\eta$ for some nontrivial ideal $J \subset \C [Z]$. Since $Z$ is homogeneous Remark \ref{adpr1} (1)-(2) implies that $Z$
has ADP and we are done.

\eproof

\providecommand{\bysame}{\leavevmode\hboxto3em{\hrulefill}\thinspace}

\end{document}